\documentclass[leqno,12pt]{article}
\usepackage{latexsym}
\usepackage{amssymb}

\usepackage{graphicx}
\usepackage{graphicx}
\usepackage{epsfig}

 0

\def \R{I \!\! R}

\newcommand{\Var}{\mbox{Var}}

\newcommand{\Cov}{\mbox{Cov}}
\newcommand{\Z}{\mathbb{Z}\mkern1mu}

\renewcommand{\theequation}{\thesection.\arabic{equation}}
\makeatletter\@addtoreset{equation}{section}\makeatother

\newcommand{\ba}{\begin{array}}
\newcommand{\ea}{\end{array}}
\newcommand{\beqohne}{\begin{eqnarray*}}
\newcommand{\eeqohne}{\end{eqnarray*}}
\newcommand{\beohne}{\begin{equation*}}
\newcommand{\eeohne}{\end{equation*}}

\renewcommand{\epsilon}{\varepsilon}

\def\3{\ss}

\newcommand{\bea}{\begin{eqnarray*}}
\newcommand{\eea}{\end{eqnarray*}}
\newcommand{\beq}{\begin{equation}}
\newcommand{\eeq}{\end{equation}}

\newcommand{\var}{ \mbox{\sl Var} \ }

\begin{document}

\title{Testing for a constant coefficient of variation in nonparametric regression}

\author{
{\small Holger Dette} \\
{\small Ruhr-Universit\"at Bochum} \\
{\small Fakult\"at f\"ur Mathematik} \\
{\small 44780 Bochum, Germany} \\
{\small e-mail: holger.dette@rub.de}\\
{\small FAX: +49 234 3214 559}\\
\and
{\small Gabriele Wieczorek} \\
{\small Ruhr-Universit\"at Bochum} \\
{\small Fakult\"at f\"ur Mathematik} \\
{\small 44780 Bochum, Germany} \\
{\small e-mail: gabriele.wieczorek@rub.de}\\
}

\maketitle

\begin{abstract}
In this paper we propose a new test for the hypothesis of a constant coefficient of variation in the common nonparametric regression model. The
test is based on an estimate of the $L^2$-distance between the square of the regression function and variance function. We prove asymptotic
normality of a standardized  estimate of this distance under the null hypothesis and fixed alternatives and the finite sample properties of
 a corresponding bootstrap test are investigated by means of a simulation study. The results are applicable to stationary processes with the common mixing conditions and are used to construct tests for ARCH assumptions in financial time series.
\end{abstract}


\section{Introduction}
\def\theequation{1.\arabic{equation}}
\setcounter{equation}{0}

We consider the common nonparametric regression model
\begin{equation} \label{1.1}
Y_i=m(X_i) + \sigma(X_i) \varepsilon_i, ~~i=1,2,\ldots,n,
\end{equation}
where $m$ denotes the regression function and $\sigma^2 $ the variance function and the random variables $\varepsilon_i$ satisfy
$E[\varepsilon_i|X_i=x]= 0$ and $E[\varepsilon_i^2|X_i=x]=1$. In many applications the variance can be assumed proportional to the squared mean
which corresponds to the assumption of a constant coefficient of variation. Typical examples include models obtained by the logarithmic
transformation from regression models with a multiplicative error structure [see Eagleson and M\"uller (1997)] or ARCH-type models [see Engle
(1982)]. Several authors have discussed the problem of estimating and testing the regression function under the restriction that $m$ and
$\sigma$ are proportional - see e.g. Mc Cullagh and Nelder (1989), who considered generalized linear models, Carroll and Ruppert (1988), who
considered a constant coefficient of variation  with a parametric model, and Eagleson and M\"uller (1997), who investigated the common
nonparametric regression model under the restriction that $m=c\sigma$ for some constant $c$.

In the present paper we will develop a formal test for the hypothesis of a constant coefficient of variation in the nonparametric regression
model (\ref{1.1}), that is
\begin{equation} \label{1.2}
H_0:  m(x) = c \sigma (x)
\end{equation}
for some positive (but unknown constant) $c$.  Besides the fact that this test can be used to check the assumptions for a statistical inference
in a nonparametric regression model with a constant coefficient of variation, it can also be used as an indicator of a multiplicative error
structure (if it is applied to the  squares of the data)
 and an exponentially distributed response $Y$ where $E[Y|X=x]=\sqrt{\Var [Y|X=x]} = m(x)$.
 In Section 2 we introduce the test statistic and indicate possible
 applications. Section 3 contains our main results in the case of an i.i.d. sample $\{X_i,Y_i \}_{i=1}^n$.  We prove asymptotic normality  of a
 standardized version of the test statistic under the null hypothesis and fixed alternatives. In Section 4 we extend these results in the case of stationary time series with the common mixing properties and discuss an application to test for an ARCH(1) model.
 The finite sample properties of a bootstrap version of the new test are investigated
 in Section 5 and some of the technical details for the proofs of our main results are presented
 in the Appendix in Section 6.

\section{Testing for a constant coefficient of variation in nonparametric regression}
\def\theequation{2.\arabic{equation}}
\setcounter{equation}{0}

Numerous authors have considered testing various hypotheses regarding the mean and the variance function in the nonparametric regression model
(\ref{1.1}) [see e.g. Dette and Munk (2003) and the references in this paper]. These hypotheses include  parametric and semi parametric
assumptions regarding the mean and variance function, but much less effort has been spent in investigating the relation between mean and
variance in the nonparametric regression model (\ref{1.1}). In the present paper we investigate the hypothesis (\ref{1.2}) of a constant
coefficient of variation using an estimate of the $L^2$-distance between the  variance and squared regression function. Typical examples
include multiplicative models of the form
$$
Y_t = m(X_t) \eta_t
$$
which can be written in the form (\ref{1.1}) with $ \sigma (\cdot )= \sqrt {\var (\eta_t)}  m (\cdot ) $ and $\varepsilon_t =
(\eta_t-1)/\sqrt{\Var(\eta_t)} $. Other examples include nonparametric ARCH models $X_t = \sqrt{m(X_{t-1}^2)} \eta_t$, for which the squared
process corresponds to a multiplicative times series model.

To be precise let
 $\{ X_i, Y_i\}^n_{i=1}$ denote a bivariate sample of
 observations from the nonparametric regression model (\ref{1.1}) with the same distribution
and let $\hat m$ and $\hat  \sigma^2 $  denote two nonparametric estimates of the regression and variance function, respectively, which will be
specified in the following section. For any positive $c$ we define the statistic $T_n(c)$ as
\begin{eqnarray}
T_n(c)&=&\frac{1}{n (n-1)} \sum_{i \ne j} K_g( X_i - X_j)  \{ c^2 Y_i^2 - (c^2 +1) \hat{m}^2(X_i) \}w(X_i)~~~~~~~~~~~~~~~~~~~~~~~   \label{2.1}\\
&~& ~~~~~~~~~~~~~~~~~~~~~~~~~~~~~~~~~\times \{ c^2 Y_j^2 - (c^2 +1) \hat{m}^2(X_j) \} w(X_j),  \nonumber
\end{eqnarray}
where $w$ denotes a weight function, $K_g(\cdot)=\frac{1}{g} K(\cdot / g)$,  $K(\cdot )$ denotes a kernel and $g$ is a bandwidth converging to
$0$ with increasing sample size.  Note that the statistic of the form (\ref{2.1}) has been considered before by Zheng (1996) for testing the
parametric form of the regression function, by Dette (2002) for testing homoscedasticity,  by Dette and von Lieres und Wilkau (2003) and Gozalo
and Linton (2000) for testing additivity in a nonparametric regression model (\ref{1.1}) with a multivariate predictor. If the estimate
$\hat{m}$ is consistent it is intuitively clear that for a large sample size
\begin{eqnarray}\nonumber
E[T_n(c)] &\approx&  E [ K_g (X_1 - X_2) \{ c^2 \sigma^2 (X_1) \varepsilon^2_1 - 2c^2 m (X_1)\sigma(X_1) \varepsilon_1 - m^2 (X_1)\} \\
\nonumber
&& ~~~~~\times \{ c^2 \sigma^2 (X_2) \varepsilon^2_2 - 2c^2 m (X_2) \sigma (X_2) \varepsilon_2 - m^2 (X_2)\} ] \nonumber \\
&\approx &
E\bigl[ f(X_i)  \{c^2 \sigma^2 (X_i)  - {m}^2(X_i) \}^2w^2(X_i) \bigr] \nonumber \\
&=& E\bigl[ \Delta_c^2(X_i) f(X_i) w^2(X_i) \bigr],  \label{2.3}
\end{eqnarray}
where $f$ denotes the density of $X$ and
\begin{equation} \label{2.4}
\Delta_c(x) = m^2(x) - c^2 \sigma^2(x).
\end{equation}
Note that $E\bigl[ \Delta_c^2(X_i) f(X_i) w^2(X_i) \bigr] =0$ if and only if the null hypothesis (\ref{1.2}) is satisfied. There exist a few
cases, where the constant $c$ in the statistic $T_n(c)$ is known. For example in  ARCH(1) models with standard normal distributed innovations
$\eta_t$ we have $X_t^2= a_0 + a_1 X_{t-1}^2+ (a_0 + a_1 X_{t-1}^2)(\eta_t^2 \nobreak-\nobreak1)$, which gives $c=1/\sqrt{2}$. However, in most
cases of practical interest the constant $c$ has to be estimated from the data. For this purpose we consider the least squares problem
\begin{equation} \label{2.4a}
\hat{c}^2 =\arg \min_{c \in \R^{>0} } \sum_{i=1}^n ( m^2(X_i) - c^2 \sigma^2(X_i))^2 w(X_i) = \frac{\sum_{i=1}^n m^2(X_i) \sigma^2(X_i)w(X_i)}{
\sum_{i=1}^n \sigma^4(X_i)w(X_i)}
\end{equation}
and estimate the unknown quantities on the right hand side. We define the residuals
\begin{equation} \label{res}
\hat{r}(X_i) = Y_i - \hat{m}(X_i),~~ (i=1,\ldots ,n)
\end{equation}
and the estimate
\begin{equation}\label{2.5}
\hat{c}^2 = \frac{(1/n)\sum_{i=1}^n \hat{m}^2(X_i) \hat{r}^2(X_i) w(X_i)}{(1/n)\sum_{i=1}^n (\hat{\sigma}^2(X_i))^2 w(X_i)},
\end{equation}
where the variance estimate $\hat \sigma$ is defined in (\ref{shat}). Note that the squared residuals $\hat{r}^2(\cdot)$ are used for estimating the variance function in the numerator of $\hat{c}^2$ in order
to avoid an additional bias caused by the use of the variance estimator $\hat{\sigma}^2(\cdot)$ [see the proof of Theorem 3.2 in the Appendix].

It is intuitively clear that the expression $\hat{c}^2$ estimates
\begin{equation} \label{2.6}
c_0^2= \frac{E[m^2(X) \sigma^2(X) w(X)]}{E[\sigma^4(X) w(X)]},
\end{equation}
which coincides with the constant $c^2$ if the null hypothesis (\ref{1.2}) is satisfied and corresponds to the best $L^2$-approximation of
$m^2$ by functions of the form $c^2 \sigma^2$, otherwise. Consequently the hypothesis of a constant coefficient of variation will be rejected
for large values of the statistic $T_n (\hat c) $.

In the following sections we specify the asymptotic properties of the statistics $T_n ( c) $, $\hat c^2  $ and $T_n (\hat c) $ if the local
linear estimate [see Fan and Gijbels (1996)] is used for estimating the mean and variance function.

\section{Main results}
\def\theequation{3.\arabic{equation}}
\setcounter{equation}{0}

In order to state our main results we have to specify nonparametric estimates of the regression and variance function and several assumptions
for the model (\ref{1.1}). We begin with the definition of the estimates. For the regression function we use the local linear estimate [see Fan
and Gijbels (1996)]
\begin{equation} \label{mhat}
\hat{m}(x)= \frac{\sum_{i=1}^n K_h(X_i -x) \left[s_{n,2}(x) - (x - X_i)s_{n,1}(x) \right] Y_i }{\sum_{i=1}^n K_h(X_i -x) \left[s_{n,2}(x) - (x
- X_i)s_{n,1}(x) \right] }
\end{equation}
where  $K_h(\cdot)=\frac{1}{h} K(\cdot / h)$, $K(\cdot )$ is a kernel, $h$ denotes a further bandwidth and
\begin{equation} \nonumber
s_{n,l}(x) = \sum_{i=1}^n K_h(X_i - x) (x - X_i)^l~~~~l=1,2.
\end{equation}
Similarly, the estimate of the variance function is obtained by replacing the observations $Y_i$ by the squared residuals $ \hat r^2 (X_i) $
defined in (\ref{res}) and is given by
\begin{equation} \label{shat}
\hat{\sigma}^2 (x)= \frac{\sum_{i=1}^n K_h(X_i -x) \left[s_{n,2}(x) - (x - X_i)s_{n,1}(x) \right] \hat r^2 (X_i) }{\sum_{i=1}^n K_h(X_i -x)
\left[s_{n,2}(x) - (x - X_i)s_{n,1}(x) \right] }.
\end{equation}
 For the sake of transparency we first assume that $\{ X_i, Y_i\}^n_{i=1}$ is a sample of independent identically distributed
observations. A corresponding result in the time series context is given in the following section. Moreover, the same bandwidths are assumed
for the calculation of the estimates of the regression and variance function for the sake of simple notation. The treatment of different
bandwidths in these estimates does not cause additional difficulties (and in the simulation study presented in Section 5 we used in fact
different bandwidths). Throughout this section we assume that the following assumptions are satisfied
\begin{itemize}
\item[(A1)] The density $f$ is twice continuously differentiable on compact sets.
\item[(A2)] The regression function $m$ is four times continuously differentiable on compact sets.
\item[(A3)] The variance function $\sigma^2$ is positive and twice continuously differentiable on compact sets.
\item[(A4)] The weight function $w$ is twice continuously differentiable and has compact support contained in $\left\{x| f(x) > 0\right\}$.
\item[(A5)] \label{kernfunktionlipschitz} The kernel
$K$ is of order $2$, and satisfies a Lipschitz condition.
\item[(A6)] \label{annahme 1} If  $n \rightarrow \infty$ the bandwidth  $g$ and $h$ satisfy
\begin{equation}\nonumber
h \sim n^{-1/5}, ~~g=o(h^2),~~ng \rightarrow \infty.
\end{equation}
\item[(A7)] The function  $m_k(x)=E[\varepsilon^k|X=x]$  is continuous for $k=3,4$ and for
$ 1 \leq k \leq 8$ uniformly bounded, that is
\begin{equation} \nonumber
E[\varepsilon^k_t|X_t=x] \leq C < \infty, ~~~~~~~~k \leq 8.
\end{equation}
\item[(A8)] The regression and variance function satisfy
\[ E[m(X)]^k < \infty~~\mbox{for}~ k=2,4, ~~\mbox{and}~~~E[\sigma^2(X)]^k < \infty~~~\mbox{for}~ k=1,2.\]
\end{itemize}
Our first result specifies the asymptotic distribution of the statistic $T_n(c)$, where the constant $c$ in the hypothesis (\ref{1.2}) is
known. Roughly speaking the statistic $T_n(c)$ is asymptotically normally distributed with different rates of convergence under the null
hypothesis and alternative. The proof is complicated and therefore deferred to the Appendix.
\bigskip

{\bf Theorem 3.1.} {\it Assume that  the assumptions (A1) - (A7) are satisfied.
\begin{itemize}
\item[(a)] Under  the null
hypothesis (\ref{1.2}) we have
\begin{equation} \nonumber
n \sqrt{g}~ T_n(c) \stackrel{\mathcal{D}}{\longrightarrow} \mathcal{N}(0,\mu_0^2),
\end{equation}
where the asymptotic variance is given by
\begin{equation}\label{1.2a}
\mu_0^2=2~  E[\{ -1 + 4 c^2 + 4 c m_3(X) + m_4(X) \}^2 m^8(X) f(X) w^4(X)]~\int K^2 (u) du.
\end{equation}
\item[(b)] Under  a fixed alternative $H_1: m \not = c \sigma $  we have
\begin{equation} \nonumber
\sqrt{n} \left( T_n(c) - E[T_n(c)] \right) \stackrel{\mathcal{D}}{\longrightarrow} \mathcal{N}(0,\mu_1^2(c)),
\end{equation}
where
\begin{equation} \nonumber
E[T_n(c)]= E[\Delta_c^2(X) f(X)w^2(X) ] + h^2 B(c) + o (h^2)
\end{equation}
with $\Delta_c$ defined in (\ref{2.4}), $\kappa_2 = 
 \int u^2 K(u) du$ and
 \begin{equation} \nonumber
B(c) = 2 (c^2 +1) ~\kappa_2 E[\Delta_c(X) m(X) m''(X) f(X) w^2(X)].
 \end{equation}
The asymptotic variance is given by
\begin{eqnarray*}
\mu_1^2(c) &=& 4 \Var(\Delta_c^2(X) f(X) w^2(X))+ 16 E[\Delta_c^2(X)m^2(X) \sigma^2(X) f^2(X) w^4(X)]\\
&~& + 4  c^4 E[ \Delta_c^2(X)\sigma^4(X) f^2(X) \{m_4(X) -1\} w^4(X)]\\
&~& - 16 c^2 E[\Delta_c^2(X) m(X) \sigma^3(X)  f^2(X) m_3(X) w^4(X)].
\end{eqnarray*}
\end{itemize}
}

\bigskip

In most applications the value $c$ in the hypothesis (\ref{1.2}) is not known and has to be estimated from the data. The following results
specify the asymptotic properties of the estimate $\hat c^2$ defined in (\ref{2.5}) and the test statistic $T_n(\hat c) $.

\bigskip

{\bf Theorem 3.2.}  {\it If the assumptions (A1) - (A8) are satisfied, then
\begin{eqnarray} \label{3.1}
 \hat{c}^2 - E [\hat{c}^2] &=& \frac{1}{n} \sum_{i=1}^n \biggl\{ \tau_1 \Bigl( m^2(X_i) \sigma^2(X_i) w(X_i) \varepsilon_i^2 - E[m^2(X) \sigma^2(X) w(X)] \Bigr)  \biggr. \nonumber \\
&~&~~~~~~~~~ + 2 \tau_1 m(X_i) \sigma^3(X_i) w(X_i) \varepsilon_i -  \tau_2   \Bigl( \sigma^4(X_i) w(X_i) - E[\sigma^4(X) w(X)] \Bigr)  \\
&~&~~~~~~~~~ - ~2 \tau_2 ~ \sigma^4(X_i) w(X_i) \{ \varepsilon_i^2 -1 \} \biggr\} + \biggl.  o_p \left( \frac{1}{\sqrt{n}}\right).\nonumber
\end{eqnarray}
Moreover,
\begin{equation} \label{3.3}
\sqrt{n} (\hat{c}^2- E[\hat{c}^2] ) \stackrel{\mathcal{D}}{\rightarrow} \mathcal{N}(0,\nu^2),
\end{equation}
where
\begin{equation} \nonumber
E[\hat{c}^2] = c_0^2+ h^2 \Gamma + o(h^2)
\end{equation}
and the constants $\Gamma, \tau_1$, $\tau_2$ and $\nu^2$ are given by
\begin{eqnarray*}
\Gamma &=& \kappa_2~ E[\sigma^2(X) \{\tau_1 m(X) m''(X)- \tau_2(\sigma^2(X))''  \}  w(X)],\\
\tau_1 &=& \frac{1}{E[\sigma^4(X) w(X)]}, \\
\tau_2&=& \frac{E[m^2(X) \sigma^2(X) w(X)]}{E^2[\sigma^4(X) w(X)]}, \\
&& \\
\nu^2 &=& \tau_1^2 \Var (m^2(X) \sigma^2(X)w(X) ) + 4 \tau_1^2 E[m^2(X) \sigma^6(X) w^2(X)] \\
&~&+ 4 \tau_1^2  E[m^3(X) \sigma^5(X) m_3(X) w^2(X)] + \tau_2^2 \Var(\sigma^4(X) w(X) )\\
&~&+ 4 \tau_2^2 E[\sigma^8(X)\{m_4(X) -1\}  w^2(X)] - 2 \tau_1 \tau_2 \Cov(m^2(X) \sigma^2(X) w(X), \sigma^4(X) w(X))\\
&~& - 4 \tau_1 \tau_2 E[m^2(X) \sigma^6(X)\{m_4(X)-1\}~  w^2(X)] - 4 \tau_1 \tau_2 E[m(X) \sigma^7(X)m_3(X) w^2(X)].
\end{eqnarray*}
}

\bigskip

We are now in a position to investigate weak convergence of the statistic $T_n(\hat c)$, where the estimate $\hat c^2$ is defined in
(\ref{2.4a}). We begin with the asymptotic distribution under the null hypothesis (\ref{1.2}). Interestingly, in this case the estimation  of
the scaling factor $c$ has no influence on the asymptotic properties of the test statistic.

\bigskip

{\bf Theorem 3.3.}  {\it Assume that  the assumptions (A1) - (A8) are satisfied. Under the null hypothesis (\ref{1.2}) we have
$$
n \sqrt{g} ~ T_n(\hat{c}) = n \sqrt{g} ~ T_n(c)  +  o_p \left( 1 \right)
 \stackrel{\mathcal{D}}{\longrightarrow} \mathcal{N}(0,\mu_0^2),
 $$
where the constant $\mu_0^2$ is defined in (\ref{1.2a}).
 }

\bigskip

Our final result in this section refers to the asymptotic properties of the statistic
 $T_n(\hat c )$ under the alternative. In this case there appears an additional term in the bias and variance of the test statistic, which is caused by
the  estimation of the scaling  factor $ c$. Recall that the constant $c^2_0$ corresponds to the best $L^2$-approximation of  $m^2$ by
functions of the form $c^2\sigma^2$.

\bigskip

{\bf Theorem 3.4.}  {\it Assume that  the assumptions (A1) - (A8) are satisfied. Under a fixed alternative
$$
\varrho =E[\Delta_{c_0}(X) \sigma^2(X) f(X) w^2(X)] > 0
$$
we have
$$
\sqrt{n} \left( T_n(\hat{c}) - E[T_n(\hat{c})]  \right) \stackrel{\mathcal{D}}{\longrightarrow} \mathcal{N}(0,\omega_1^2),
$$
where
$$
E[T_n(\hat{c})]= E[\Delta_{c_0}^2(X) f(X)w^2(X)] + h^2( B(c_0) - 2 \varrho \Gamma )+ o(h^2),
$$
and $B(c_0)$ is a term in the bias of the statistic $T_n(c_0)$. The asymptotic variance $\omega_1^2$ is given by
$$
\omega_1^2 = \mu_1^2(c_0) + 4 \varrho^2 \nu^2 - 4 \varrho~ \upsilon^2(c_0),
$$
where  $\mu_1^2(c_0)$ is defined in Theorem 3.1(b), $\nu^2$ corresponds to the asymptotic variance of $\hat{c}^2$ in Theorem 3.2 and
\begin{eqnarray*}
\upsilon^2(c_0) &=& 2 \tau_1 E[ \Delta_{c_0}(X) (m^2(X) - c_0^2 \sigma^2(X) m_4(X)) m^2(X) \sigma^2(X) f(X) w^3(X)  )\\
&~& - 2 \tau_1 E[\Delta_{c_0}^2(X) f(X) w^2(X)] E[m^2(X) \sigma^2(X) w(X)]\\
&~&- 4  c_0^2 \tau_1 E[ \Delta_{c_0}(X) m(X) \sigma^5(X)  f(X) m_3(X)w^3(X)  ]\\
&~& - 2 \tau_2 \Cov (\Delta_{c_0}^2(X)f(X) w^2(X), \sigma^4(X) w(X))\\
&~& + 4  c_0^2  \tau_2 E [ \Delta_{c_0}(X)  \sigma^6(X) f(X)\{m_4(X) -1\} w^3(X) ]\\
&~&+ 4 \tau_1   E [\Delta_{c_0}(X) m^3(X) \sigma^3(X) f(X)m_3(X) w^3(X) ] \\
&~& + 8  \tau_1 E[ \Delta_{c_0}(X) m^2(X) \sigma^4(X) f(X) w^3(X)] \\
&~& - 8 \tau_2  E[\Delta_{c_0}(X) m(X) \sigma^5(X) f(X)m_3(X) w^3(X)].
\end{eqnarray*}
}

{\bf Remark 3.5.} The term $\upsilon^2(c_0) $ corresponds to the asymptotic covariance between the statistic $T_n(c_0)$  and the estimate $\hat
c^2$ of $c^2_0$.

\section{Further discussion}
\def\theequation{4.\arabic{equation}}
\setcounter{equation}{0}

\subsection{Asymptotic results for absolutely regular processes}

The general nonparametric framework includes time series models. Typical examples are multiplicative models $Z_t =\sigma_t \eta_t$, where
$\sigma_t$ is a positive function of the past $\{Z_{t-i}: i \geq 1\}$ and possibly of the past volatility $\{\sigma_{t-i}: i \geq 1\}$. For
instance, defining $\sigma_t$ by $\sqrt{ \vartheta_0 + \vartheta_1 Z_{t-1}^2}$ for $\vartheta_i \geq 0$ we achieve the linear ARCH(1) model.
Therefore our test can also be used as a preliminary step to identify certain time series. For this purpose it is necessary to extend the
asymptotic results under a more general setup which includes both time series data and i.i.d. observations as special cases. For this purpose
we need the following assumptions for some fixed $\varepsilon \in (0,1/2)$ and $\xi > 2$.

\begin{itemize}
\item[(M1)] The process ${(X_i, Y_i)}$ is absolutely regular, i.e.
\[ \beta(k) = \sup_{s \in \Z} E[\sup \{ |P(A| \mathcal{F}_{-\infty}^s) - P(A) | A \in \mathcal{F}_{s+k}^{\infty}\}] \rightarrow 0, ~~\mbox{as}~~~ k \rightarrow 0,
\]
where $\mathcal{F}_s^t$ is the $\sigma$-algebra generated by $\{(X_l, Y_l): s \leq l \leq t\}$. Further,
\[\sum_{j=1}^{\infty} j^2 \beta^{\frac{\epsilon}{1 + \epsilon}}(j) < \infty. \]
\item[(M2)] The innovations $\varepsilon_t$ in the model (\ref{1.1}) satisfy
\[ E[\varepsilon_t|X_t , \mathcal{F}_{- \infty}^{t-1}(X,Y)] = E[\varepsilon_t|X_t]=0, \]
and \[ \Var(Y_t| X_t, \mathcal{F}_{- \infty}^{t-1}(X,Y) ) = \sigma^2(x).\] Further, $E|\varepsilon_t|^k < \infty$ to the order  $k \leq 48 \xi
(1+ \varepsilon)$.
\item[(M3)] The regression function $m(\cdot)$ satisfies
\begin{eqnarray*}
&~&E|m(X)|^k < \infty ~~\mbox{for $k \leq 4 (1+ \varepsilon)$}~~~~\mbox{and}~~~~~~~~~~~~~~~~~~~~~~~~~~\\
&~&E|m''(X)|^k < \infty ~~~\mbox{for $k \leq 20 \xi (1+ \varepsilon)$},
\end{eqnarray*}
whereas the variance function $\sigma^2(\cdot)$ fulfills
\[ E|\sigma^2(X)|^k < \infty ~~~~~~~\mbox{for $k \leq 12 \xi (1+ \varepsilon).~~~~$}~~~~~~~~~~~~~~~\]
\end{itemize}

Note that  assumption (M3) contains  assumption (A8) which is therefore omitted. Under the assumptions (A1) - (A7) together with (M1) - (M3)
the asymptotic results for $T_n(c), \hat{c}^2$ and $T_n(\hat{c})$ can be established for strictly stationary, $\beta$-mixing processes
$\{X_i,Y_i\}_{i \in \Z}$. The proof of the following results is obtained from the proof of the statements presented in Section 3 for the
independent case using similar arguments as given by Dette and Spreckelsen (2004), where the authors investigate the asymptotic distribution of
goodness-of-fit tests of linearity for absolutely regular processes. For the sake of brevity the details are omitted and we refer the
interested reader to the PhD thesis of Wieczorek (2007). Moreover,  we only state the results for the statistic $T_n(\hat{c})$.
Note that under the null hypothesis the asymptotic distribution of $T_n(\hat{c})$ under mixing assumptions coincides with the distribution for the i.i.d. case.\\

{\bf Theorem 4.1} {\it Assume that the assumptions (A1) - (A7) and (M1) - (M3) are satisfied. Under the null hypothesis, we have
\[n \sqrt{g} ~T_n(\hat{c})  \stackrel{\mathcal{D}}{\longrightarrow} \mathcal{N}(0,\mu_0^2), \]
where $\mu_0^2$ is the asymptotic variance of $T_n(c)$ defined in (\ref{1.2a}).
}\\

Our final theoretical result states the asymptotic properties of the statistic $T_n(\hat{c})$ under fixed alternatives. Note that in this case
the variance of the limit distribution contains the variance of the limit distribution for the i.i.d. case as well as additional covariances.
For a precise statement of the result we introduce the notation $E^{\otimes}$,
which denotes the expectation with respect to the product measure. \\

{\bf Theorem 4.2.} {\it If the assumptions (A1) - (A7) and (M1) - (M3) are satisfied, then under a fixed alternative $\varrho > 0$ we have
\[ \sqrt{n} \left( T_n(\hat{c}) - E^{\otimes}[T_n(\hat{c})] \right) \stackrel{\mathcal{D}}{\longrightarrow} \mathcal{N}(0,\tilde{\omega}_1^2). \]
In particular,
\[ \left|E[T_n(\hat{c})]- E^{\otimes}[T_n(\hat{c})]\right| = o\left( \frac{1}{\sqrt{n}}\right), \]
where the mean $E^{\otimes}[T_n(\hat{c})]$ and the constant $\varrho$ are given in Theorem 3.4. The asymptotic variance is given by
\begin{equation} \label{asymvariancedependent}
 \tilde{\omega}_1^2 = \tilde{\mu}_1^2(c_0) + 4 \varrho^2 \tilde{\nu}^2 - 4 \varrho~ \tilde{\upsilon}^2(c_0),
 \end{equation}
where $\tilde{\mu}_1^2(c_0)$ denotes the asymptotic variance of $T_n(c_0)$ defined by
\begin{eqnarray*}
\tilde{\mu}_1^2(c_0) &=& \mu_1^2(c_0) + 8 \sum_{t=1}^{\infty} \Cov \bigl( \Delta_{c_0} (X_1)  [\Delta_{c_0}(X_1, \varepsilon_1) + 2 m(X_1) \sigma(X_1) \varepsilon_1]f(X_1) w^2(X_1),\bigr.\\
&~&~~~~~~~~~~~~~~~~~~~~~~~~~\bigl. \Delta_{c_0}^2(X_{1+t}) f(X_{1+t})w^2(X_{1+t})\bigr),
\end{eqnarray*}
$\mu_1^2(c_0)$ is defined in Theorem 3.1(b). The term $\tilde{\nu}^2$ in (\ref{asymvariancedependent}) corresponds to the asymptotic variance
of the estimate $\hat{c}^2$ given by
\begin{eqnarray*}
\tilde{\nu}^2&=& \nu^2 + 2 \sum_{t=1}^{\infty} \Cov\Bigl( 2 \tau_1 m(X_1) \sigma^3(X_1) w(X_1) \varepsilon_1-  \tau_2   \sigma^4(X_1) w(X_1) \{ 2 \varepsilon_1^2 -1\}   \Bigr.\\
&~&~~~~~~~~~~~~~~~~~~~~\Bigl.+ \tau_1 m^2(X_1) \sigma^2(X_1) w(X_1) \varepsilon_1^2, \tau_1 m^2(X_{1+t}) \sigma^2(X_{1+t}) w(X_{1+t}) -  \tau_2
\sigma^4(X_1) w(X_1) \Bigr),
\end{eqnarray*}
where $\nu^2$ is given in Theorem 3.2 and $\tilde{\upsilon}^2(c_0)$ corresponds to the asymptotic covariance between $T_n(c_0)$ and $\hat{c}^2$
defined by
\begin{eqnarray*}
\tilde{\upsilon}^2 (c_0)&=&\upsilon^2(c_0) + 2 \sum_{t=1}^{\infty}  \Cov \left(\Delta_{c_0}(X_1)f(X_1) \left\{ \Delta_{c_0}(V_1)+ 2 m(X_1) \sigma(X_1) \varepsilon_1  \right\}w^2(X_1), \right.\\
&~& ~~~~~~~~~~~~~~~~\left. \left\{\tau_1 m^2(X_{1+t}) \sigma^2(X_{1+t})  - \tau_2 \sigma^4(X_{1+t})\right\} w(X_{1+t}) \right) \\
&~&~~~~~~~~+ 2 \sum_{t=1}^{\infty} \Cov \left(\Delta_{c_0}^2(X_{1+t}) f(X_{1+t}) w^2(X_{1+t}), \right.\\
&~& ~~~~~~~~~~~~~~~~ \left. \left\{\tau_1  m^2(X_1) \sigma^2(X_1)  \varepsilon_1^2 + 2\tau_1  m(X_1) \sigma^3(X_1)\varepsilon_1- \tau_2
\sigma^4(X_1) ( 2 \varepsilon_1^2 -1)  \right\} w(X_1) \right),
\end{eqnarray*}
and $\upsilon^2(c_0)$ is defined in Theorem 3.4.
}\\

{\bf Remark 4.3.}  It is worthwhile to mention that in the case where the stationary process is absolutely regular with a geometric
rate, i.e. $\beta (j) = O( \rho ^j)$ for some $\rho \in (0,1)$, the asymptotic covariance of the test statistic given in
Theorem 4.2 coincides with the asymptotic covariance given in Theorem 3.4 for the independent case, that is:
$$
\tilde \mu_1^2(c_0) =\mu_1^2(c_0)~,~~\tilde \nu^2 =\nu^2~,~~\tilde \nu^2(c_0) =\nu^2(c_0)~.
$$
\bigskip

{\bf Remark 4.4.} The moment assumption (M3) is quite restrictive and limits the applicability of the  test to many interesting time series
models such as ARCH or GARCH models. One possible way to circumvent assumption (M3) is the introduction of an additional weight function in the
estimates. As a consequence a slight modification of the estimates can be arranged in our testing procedure eliminating assumption (M3). The
details can be found in Wieczorek (2007), and only the modification is mentioned for the sake of brevity. We introduce in a first step an
additional weight function $w^{\ast}$, satisfying

\begin{itemize}
\item[(A9)] $w^{\ast}$ is twice continuously differentiable and has compact support contained in $\left\{x| w(x) > 0\right\}$.
\end{itemize}

Next, we propose a modified estimate of the regression function given by $\hat{m}^{\ast}(x)= \hat{a}$, where
\begin{equation} \label{leastsquaresw}
 (\hat{a}, \hat{b})= \mbox{arg}\min_{a,b} \sum_{i=1}^n \left\{ Y_i - a - b(X_i -x) \right\}^2 w(X_i) K\left( \frac{X_i-x}{h}\right)
\end{equation}
is the local linear estimate (additionally weighted by $w$) of the regression function and its derivative. Note that the modified local linear
regression estimator $\hat{m}^{\ast}$ differs from the local linear estimate $\hat{m}$ in (\ref{mhat}) by the introduction of the weight
function $w$ in (\ref{leastsquaresw}). Similarly, we propose $(\sigma^{ 2})^\ast(x)=\hat{\alpha}$ as the modified estimate of the variance
function, where
\[ (\hat{\alpha}, \hat{\beta})= \mbox{arg}\min_{\alpha,\beta} \sum_{i=1}^n \left\{ (\hat{r}^{\ast})^2(X_i) - a - b(X_i -x) \right\}^2 w^{\ast}(X_i) K\left( \frac{X_i-x}{h}\right) \]
is the local linear estimate (weighted by the second weight function $w^{\ast}$) based on the nonparametric residuals $\hat{r}^{\ast}(X_i)$
defined by
\[ \hat{r}^{\ast}(X_i) = Y_i -\hat{m}^{\ast}(X_i).  \]

Based on the modified estimates of the regression function and the variance function the new test statistic is defined by
\begin{eqnarray}
T_n^{\ast}(c)&=&\frac{1}{n (n-1)} \sum_{i \ne j} K_g( X_i - X_j)  \left\{ c^2 Y_i^2 - (c^2 +1) (\hat{m}^{\ast})^2(X_i) \right\}w^{\ast}(X_i)~~~~~~~~~~~~~~~~~~~~~~~ \nonumber \\
&~& ~~~~~~~~~~~~~~~~~~~~~~~~~~~~~~~~~~~~\times \left\{ c^2 Y_j^2 - (c^2 +1) (\hat{m}^{\ast})^2(X_j) \right\} w^{\ast}(X_j).
\label{teststatistikast}
\end{eqnarray}
In addition, we consider the modified least squares problem
\[ (\hat{c}^2)^{\ast} =\arg \min_{c \in \R^{>0} } \sum_{i=1}^n ( m^2(X_i) - c^2 \sigma^2(X_i))^2 (w^{\ast})^3(X_i).  \]
Therefore, we define the estimate of $c^2$ by
\[ (\hat{c}^2)^{\ast} = \frac{(1/n)\sum_{i=1}^n (\hat{m}^{\ast})^2(X_i) (\hat{r}^{\ast})^2(X_i) (w^{\ast})^3(X_i)}{(1/n)\sum_{i=1}^n ((\hat{\sigma}^2)^{\ast}(X_i))^2 (w^{\ast})^3(X_i)}.\]
As an immediate consequence of the modified definitions the asymptotic results in Theorem 4.1 and 4.2 can also be established for the test
statistic $T_n^{\ast}(\hat{c}^{\ast})$. No additional assumptions are needed, in particular, the introduction of the weight functions in the
estimators avoids the assumption (M3) about the boundedness of the moments of the regression and variance function [see Wieczorek (2007) for
details]. This modification makes the test applicable to financial time series, as demonstrated in the following section.

\bigskip

{\bf Remark 4.5.} An alternative test statistic, which is  related closer to the estimation of the coefficient of variation $m(x)/\sigma(x)$ was
proposed by the referee. This approach is similar 
to the method discussed in Dette and Hetzler (2008) and is currently being investigated by
the authors. To be precise, define $U_i = Y_i /\hat \sigma (X_i)$ $(i=1,\dots,n)$, $\overline U = \frac {1}{n} \sum^n_{j=1} U_j$ and consider
the statistic
$$
H_n (x) = \int^x_{- \infty} \frac {1}{n}  \sum^n_{j=1} w (X_j) K_h (X_j - t) (U_j - \overline U) dt .
$$
A heuristic computation shows that
$$
E[H_n (x)] \approx \frac {n-1}{n} \int^x_{- \infty} w(t) f(t) \left \{ \frac {m(t)}{\sigma (t)} - \int \frac {m(s)}{\sigma (s)} f(s) ds \right
\} dt
$$
which vanishes a.e.\ if and only if the null hypothesis of a constant coefficient of variation is satisfied. Consequently, a test for this
hypothesis could be based on a Kolmogorov-Smirnov- or Cr\'{a}mer-von Mises statistic of the process $H_n$. For a similar approach in the
context of testing for a parametric form of the variance function we refer to Dette and Hetzler (2008).

\subsection{Example: Application to financial time series}

The hypothesis of the proportionality of the regression function $m$ and the volatility function $\sigma$ can also be used to test for a
multiplicative model structure. In particular, the proposed test can be viewed as a preliminary step in time series analysis before applying
other procedures such as specific testing procedures for ARCH or GARCH models. One important criterion in order to establish all asymptotic
results in such a context is assumption (M3). There the existence of bounds for the absolute moments of the regression function $m$, its second
derivative $m''$ and the variance function $\sigma^2$ is required. But often financial time series do not satisfy this assumption. For
instance, consider the linear ARCH(1) model
\[ Z_t = \sqrt{\vartheta_0 + \vartheta_1 Z_{t-1}^2} \eta_t \]
for some constants $\vartheta_0, \vartheta_1 \geq 0, \vartheta_1 < 1$, where $\eta_t$ has mean 0 and variance 1 and is independent of $Z_{t-1}$
for all $t$. The squared ARCH(1) process can be written as
\begin{equation} \label{archquad}
Z_t^2 = ( \vartheta_0 + \vartheta_1 Z_{t-1}^2) + (\vartheta_0 + \vartheta_1 Z_{t-1}^2) \varepsilon_t,
 \end{equation}
 where $\varepsilon_t = \eta_t^2 -1$. Clearly,  model (\ref{archquad}) can be identified as a particular case of the general nonparametric
 regression model (\ref{1.1}) by taking $Y_t= Z_t^2$, $X_t=Z_{t-1}^2$, $m(X_t)= \vartheta_0 + \vartheta_1 X_t$
 and $\sigma(X_t) = c^{-1} (\vartheta_0 + \vartheta_1 X_t)$. The scaling factor $c$ is given by $c^2 = (E[\eta^4] -1 )^{-1}$ and
 depends only on the error distribution.\\
For the ARCH(1) process the assumption (M3) can therefore be formulated in terms of the boundedness of absolute moments of $Z_t$. So it is
important to know whether the stationary solution $Z_t$ has moments of higher orders to apply the test. For stationary ARCH(p) processes with a
symmetric error distribution, a necessary and sufficient condition for the existence of such moments  has been given by Milh\o j (1985). In
particular, let $m > 0$, then the $m$'th moment of an ARCH(1) model exists if and only if $E[\vartheta_1 \eta_0^2]^m < 1$.
As an immediate consequence, one sees that in many cases ARCH processes do not have finite moments of higher orders. \\
In such cases we refer to Remark 4.4. In order to circumvent the assumption of existing high-order moments of $Z_t$ we apply the (slightly)
modified testing procedure.  In particular, the identification of the regression function $m$ and the variance function $\sigma^2$ provides the
assumptions (A2) and (A3) to be satisfied. Furthermore, from
\[ E[\varepsilon^k |X = x] = c^k E[(\eta^2 - 1)^k| Z_{t-1}^2 ]= c^k E[(\eta^2 - 1)^k ]\]
it follows that $(\varepsilon_t)$ fulfills (A7) and (M2) if the innovations $\eta_t$ satisfy certain moment conditions. If the assumptions
(A1), (A4) - (A7), (A9) are satisfied and the ARCH(1) process $(Z_t)$ fulfills the assumption (M1) the asymptotic normality under the null
hypothesis of a multiplicative model of the corresponding test statistic $T_n^{\ast}(\hat{c}^{\ast})$ can be established, that is
\[ n \sqrt{g} \left(T_n^{\ast}(\hat{c}^{\ast}\right) - E [T_n^{\ast}(\hat{c}^{\ast})] ) \stackrel{\mathcal{D}}{\rightarrow} \mathcal{N}(0, (\mu_0^2)^{\ast}), \]
where the asymptotic variance $(\mu_0^2)^{\ast}$ is given by
\[ (\mu_0^2)^{\ast} = 1152 \int K^2(u) du ~~\sum_{k=0}^8 {8 \choose k}  \vartheta_0^{n-k} \vartheta_1^k E[Z^{2k} f(Z^2) (w^{\ast})^4(Z^2)].
\]

\section{Finite sample properties}
\def\theequation{5.\arabic{equation}}
\setcounter{equation}{0}

In order to study the finite sample properties of the new test we have conducted a small simulation study. Because it is well known that the
approximation of the nominal level by the normal distribution provided by Theorem 3.3 is not very accurate for moderate sample sizes, we do not
recommend to estimate the asymptotic variance and bias and to compare the standardized statistic with the quantiles of a normal distribution.
Instead, we propose to use resampling methods. As an example we
 have
implemented a smooth bootstrap procedure to obtain the critical values. For this purpose we estimate the regression and variance function by
the local linear estimates defined in (\ref{mhat}) and (\ref{shat}), respectively, and consider the standardized residuals
\begin{equation} \label{s1}
\eta_i = \frac {Y_i - \hat m (X_i)}{\hat{\sigma} (X_i)} \quad i = 1, \dots, n \: ,
\end{equation}
which are normalized to have mean 0 and variance 1, that is
\begin{equation} \label{s2}
\hat \varepsilon_i = \frac {\eta_i - \overline \eta}{\sqrt{\frac {1}{n-1} \sum^n_{i=1} (\eta_i - \overline \eta)^2}} \quad i = 1, \dots, n \: .
\end{equation}
The bootstrap errors are then defined as
\begin{equation} \label{s3}
\varepsilon^\ast_i = \tilde \varepsilon^\ast_i + \mbox{v} N_i,
\end{equation}
where $\tilde \varepsilon^\ast_i, \dots, \tilde \varepsilon^\ast_n$ are drawn randomly with replacement from the empirical distribution of the
standardized residuals $\hat \varepsilon_1, \dots, \hat \varepsilon_n$ and $N_1, \dots, N_n$ are i.i.d\ standard normal distributed random
variables independent of the sample ${\cal Y}_n = \{(X_1, Y_1), \dots,$ $(X_n, Y_n)\}$ and v=v$_n$ is a smoothing parameter converging to 0
with increasing sample size. In the next step bootstrap data is generated according to the model
\begin{equation} \label{s4}
Y^\ast_i = \hat c \ \hat \sigma (X_i) + \hat \sigma (X_i) \varepsilon^\ast_i \quad i=1,\dots,n \: ,
\end{equation}
 where $\hat c$ is the least squares estimate (\ref{2.5}) obtained from the data corresponding to the range
 $[X_{( \lfloor 0.05n \rfloor ) },  X_{( \lfloor 0.95n \rfloor ) }]$ of the predictors,
where $X_{(1)} \leq \dots \leq X_{(n)}$ denotes the order statistic of $X_1,\dots, X_n$.
 The test statistic $T^\ast_n$ is calculated from the bootstrap data $(X_1, Y^\ast_1)$, $\dots, (X_n, Y^\ast_n)$. If B bootstrap replications have been performed, the null hypothesis (\ref{1.2}) is rejected if
\begin{equation} \label{s5}
T_n > T_n^{\ast (\lfloor B(1- \alpha) \rfloor)} \: ,
\end{equation}
where $T_n^{\ast (1)} < \quad < T_n^{\ast (B)}$ denote the order statistic of the bootstrap sample. For the size of the bootstrap replications
we chose B = 100, while 1000 simulation runs are performed for the calculation of the empirical level of this test. The sample sizes are given
by $n = 50, 100, 200$ and the smoothing parameters in the test statistic and the bootstrap procedure are chosen by $g = n^{-1/2}$ and v = 0.1,
respectively. The bandwidths for the estimation of the variance and regression function are chosen separately by least squares cross
validation.

Our first example considers the model
\begin{equation} \label{s6}
m(x) = c (1 + 0.1 x ) \: ; \quad \sigma (x) = (1 + 0.1 x),
\end{equation}
where c = 0.5, 1, 1.5. The predictors $X_1, \dots, X_n$ are independent identically distributed following a uniform distribution on the
interval [0,1], while the errors $\varepsilon_1, \dots, \varepsilon_n$ have a standard normal distribution. The first part of Table 1 shows the
approximation of the nominal level, which is rather accurate for sample sizes larger than $n = 100$.
For the smaller sample sizes $n=50$ and larger values of $c$ the approximation is not very accurate.
In a second step we study the power of the
test and consider the models
\begin{equation} \label{s7}
m(x) = c (1 + 0.1 x ) \: ; \quad \sigma (x) = (1 + 0.1 x + \sqrt{x}) \: ,
\end{equation}
\begin{equation} \label{s8}
m(x) = c (1 + 0.1 x ) \: ; \quad \sigma (x) = (1 + 0.1 x + 2\sqrt{x}) \: .
\end{equation}
The corresponding results are depicted in the lower part of Table 1. For the model (\ref{s7}) we observe a moderate increase in power, which
corresponds to intuition. Because the predictor varies in the interval [0,1], the deviation from a multiplicative structure is extremely small
for model (\ref{s7}). On the other hand, the alternative model (\ref{s8}) is detected with larger power, which is also reflected by rather high
simulated rejection probabilities.

\begin{table}[h] \begin{small}
\begin{center}
\begin{tabular}{|c|c|r|r|r|r|r|r|r|r|r|r|r|r|}
\hline
\multicolumn{2}{|c|}{$n$} &\multicolumn{4}{|c|}{$50$}&\multicolumn{4}{|c|}{$100$}  &\multicolumn{4}{|c|}{$200$}\\
\hline
&$c$ $\: \backslash \: \alpha$ &2.5 \% &5 \% &10 \%&20 \% &2.5 \% &5 \% &10 \% &20 \% &2.5 \% & 5 \%&10 \% &20 \%\\
\hline &0.5&.028&.056&.104&.229&.035&.063&.108&.211&.032&.052&.099&.207\\ \cline{2-14} (\ref{s6}) &1.0 &.042 &.051 &.105 &.210 &.031 &.049
&.099 & .194 &.033&.054&.104&.204\\ \cline{2-14}
&1.5&.051&.069&.123&.241&.041&.061&.109&.202&.038&.055&.096&.205\\
\hline &0.5& .043&.067&.144&.267 &.051 &.074&.148&.269& .122 &.187&.281&.442 \\ \cline{2-14} (\ref{s7}) &1.0&.041&.069&.136&.264& .044&.094
&.168 &.281 &.105&.144&.233& .368\\ \cline{2-14} &1.5&.064&.091&.141&.247&.077&.109&.166&.273 &.114 & .152 &.208& .341 \\  \hline
&0.5&.049&.092&.167&.283&.066&.097&.213&.368 &.097 &.155& .261& .421\\ \cline{2-14} (\ref{s8})
&1.0&.073&.122&.215&.362&.092&.156&.297&.464&.177 & .266&.388&.554 \\ \cline{2-14} &1.5&.063&.107 &.203&.378& .123  & .187 &.281&.442 & .233&
.316& .424& .559 \\  \hline
\end{tabular}
\label{tab1} \caption{{\it Simulated rejection probabilities of the bootstrap test (\ref{s5}), for three nonparametric regression models, where
the first line corresponds to a multiplicative model.}}
\end{center} \end{small}
\end{table}

\bigskip
Our second  example investigates the performance of the bootstrap test in the context of stationary time series. To this end we consider two
models corresponding to the null hypothesis, that is
\begin{equation} \label{sta1}
X_t = ( 1 + 0.1 \ X_{t-1}) + (1 + 0.1 \ X_{t-1}) \varepsilon_t
\end{equation}
\begin{equation} \label{sta2}
X_t = \sin( 1 + 0.5 \ X_{t-1}) + \sin(1 + 0.5 \ X_{t-1}) \varepsilon_t
\end{equation}
and two models corresponding to the alternatives of no multiplicative structure, i.e.
\begin{equation} \label{sta3}
X_t = ( 1 + 0.1 \ X_{t-1}) + 0.5 \sqrt{| X_{t-1}|} \ \varepsilon_t
\end{equation}
\begin{equation} \label{sta4}
X_t = \sin( 1 + 0.5 \ X_{t-1}) + \cos (1 + 0.5 \ X_{t_1}) \varepsilon_t
\end{equation}
where the innovations are again independent standard normal distributed. The corresponding results are displayed in Table 2. We observe a
reasonable approximation of the nominal level for the two models corresponding to the null hypothesis. On the other hand, the two
 alternatives in (\ref{sta3})
and (\ref{sta4}) are detected with reasonable power.

\bigskip

\begin{table}[h] \label{tabneu}
\begin{center}
\begin{tabular}{|c|r|r|r|r|r|r|r|r|r|r|r|r|r|}
\hline
$n$&\multicolumn{4}{|c|}{$50$} &\multicolumn{4}{|c|}{$100$} &\multicolumn{4}{|c|}{$200$}\\
\hline
$\alpha$&2.5 \% & 5 \% & 10 \% & 20 \% & 2.5 \% & 5 \% & 10 \% &20 \% & 2.5 \% &5 \% &10 \%&20 \% \\
\hline
(\ref{sta1})& .029&  .047& .097& .217& .023& .048& .089& .187& .024& .048& .097& .189\\
\hline
(\ref{sta2}) & .038 & .057 & .109 & .201 & .035 & .054 & .092 & .191 & .036& .057& .109& .205\\
\hline
(\ref{sta3})& .053& .077& .161& .295& .074& .092& .182& .314& .113& .156& .237& .395\\
\hline
(\ref{sta4})& .084& .117& .189& .299& .097& .133& .212& .321& .129& .176& .289& .417\\
\hline
\end{tabular}
\caption{{\it Simulated rejection probabilities of the bootstrap test (\ref{s5}) for four nonparametric autoregressive time series models. The
models (\ref{sta1}) and (\ref{sta2}) correspond to the null hypothesis of a multiplicative model, while models (\ref{sta3}) and (\ref{sta4})
represent the alternative.}}
\end{center}
\end{table}

\section{Appendix:  proofs}
\def\theequation{6.\arabic{equation}}
\setcounter{equation}{0}

\subsection{Proof of Theorem 3.1.}

A straightforward calculation gives the decomposition
\begin{equation} \label{5.1}
T_n(c) = (c^2 + 1)^2 T_{1n} - 2(c^2 +1)  \{  2c^2T_{2n} -  T_{3n}(c) \} + T_{4n}(c) - 4 c^2  \{ T_{5n}(c) - c^2 T_{6n}\} ,
\end{equation}
with
\begin{eqnarray*}
T_{1n} &=& \frac{1}{n(n-1)} \sum_{i \ne j}  K_g( X_i - X_j)  \delta_n(X_i)w(X_i)  \delta_n(X_j)w(X_j), \\
T_{2n} &=& \frac{1}{n(n-1)} \sum_{i \ne j}  K_g( X_i - X_j)  \delta_n(X_i) w(X_i)  m(X_j) \sigma(X_j) w(X_j) \varepsilon_j, \\
T_{6n} &=& \frac{1}{n(n-1)} \sum_{i \ne j}  K_g( X_i - X_j)  m(X_i) \sigma(X_i)w(X_i) \varepsilon_i  m(X_j) \sigma(X_j) w(X_j)\varepsilon_j, \\
T_{3n}(c) &=& \frac{1}{n(n-1)} \sum_{i \ne j}  K_g( X_i - X_j)  \delta_n(X_i)w(X_i)  \Delta_c(X_j, \varepsilon_j) w(X_j), \\
T_{4n}(c) &=& \frac{1}{n(n-1)} \sum_{i \ne j}  K_g( X_i - X_j)  \Delta_c(X_i, \varepsilon_i) w(X_i) \Delta_c(X_j, \varepsilon_j)w(X_j), \\
T_{5n}(c) &=& \frac{1}{n(n-1)} \sum_{i \ne j}  K_g( X_i - X_j) \Delta_c(X_i, \varepsilon_i)w(X_i)  m(X_j) \sigma(X_j) w(X_j)\varepsilon_j,
\end{eqnarray*}
where  we have used the notation
\begin{eqnarray}
\Delta_c(X_i,\varepsilon_i)  &=&m^2(X_i) - c^2 \sigma^2(X_i) \varepsilon^2_i   \label{5.2} \\
\delta_n(X_i) &=& \hat{m}^2(X_i)- m^2(X_i). \label{5.3}
\end{eqnarray}
At the end of the proof  we will show that  the terms $T_{1n}$ and $T_{2n}$ are asymptotically negligible under the null hypothesis and under
fixed alternatives, that is
\begin{equation} \label{5.4}
n \sqrt{g}~ T_{jn} \stackrel{p}{\longrightarrow} 0, ~~~~~~j=1,2.
\end{equation}
We now have to distinguish the case of the null hypothesis and alternative.

\bigskip

{\bf Proof of Theorem 3.1(a).}  Note that the statistic $T_{3n}(c)$ reduces under the null hypothesis to
\begin{eqnarray} \label{5.4a}
T_{3n}(c) &\stackrel{H_0}{=}& \frac{1}{n(n-1)} \sum_{i \ne j}  K_g( X_i - X_j ) \delta_n(X_i) w(X_i) m^2(X_j) w(X_j)\{ 1 - \varepsilon_j^2\}. \\
& \stackrel{H_0}{=}& 2 \tilde{T}_{3n}^{(1)} + \tilde{T}_{3n}^{(2)} \nonumber
\end{eqnarray}
with
\begin{eqnarray*}
\tilde{T}_{3n}^{(1)} &=& \frac{1}{n (n-1)} \sum_{i \ne j} K_g(X_i - X_j) m(X_i) \tilde{\delta}_n(X_i) w(X_i) m^2(X_j) w(X_j) \{ 1 - \varepsilon_j^2\}, \\
\tilde{T}_{3n}^{(2)}&=& \frac{1}{n(n-1)} \sum_{i \ne j} K_g(X_i - X_j) \tilde{\delta}_n^2(X_i) w(X_i) m^2(X_j) w(X_j) \{ 1- \varepsilon_j^2 \},
\end{eqnarray*}
where  we use the notation
\begin{equation} \label{5.5aa}
\tilde{\delta}_n(x) = \hat{m}(x) - m(x).
\end{equation}
An application of Lemma 2 in Yao and Tong (2000) and a straightforward calculation [similar to the proof of the estimate (\ref{5.4}) given at
the end of this section] yields
\begin{equation} \nonumber
\tilde{T}_{3n}^{(1)} = O_p\left( \frac{1}{nh} \right)~,~~\tilde{T}_{3n}^{(2)} = O_p \left( \frac{h}{n} \right).
\end{equation}
Observing the conditions on the bandwidths we therefore have
\begin{equation} \label{5.4c}
T_{3n}(c) \stackrel{H_0}{=}o_p\left( \frac{1}{n \sqrt{g}} \right).
\end{equation}
We note again that this estimate holds only under the null hypothesis.

We now obtain from (\ref{5.4}) and (\ref{5.4c}) under the null hypothesis $H_0:$ $m(\cdot) = c~\sigma(\cdot)$ that
\begin{equation} \label{5.4d}
T_n(c)\stackrel{H_0}{=} T_{4n}(c) - 4 c^2 T_{5n}(c) + 4 c^4 T_{6n} + o_p \left( \frac{1}{n \sqrt{g}}\right),
\end{equation}
and the assertion follows if the weak convergence can be established for the statistic
\begin{eqnarray}
\tilde{T}_n(c) &=& T_{4n}(c) - 4 c^2 T_{5n}(c) + 4 c^4 T_{6n} =  \frac{2}{n(n-1)} \sum_{i < j} h_n(V_i,V_j) ,\label{5.4e}
\end{eqnarray}
where $V_i=(X_i, \varepsilon_i)$,
\begin{eqnarray*}
h_{ij} &=& h_n(V_i,V_j) =  K_g(X_i - X_j)  \pi(V_i) w(X_i)  \pi(V_j) w(X_j), \\
\pi(V_i) &=& m^2(X_i)  \{ 1- \varepsilon_i^2 - 2 c \varepsilon_i\}.
\end{eqnarray*}
For a proof of weak convergence of $\tilde{T}_n(c)$ we note that the kernel $h_n(\cdot,\cdot)$  is degenerate, i.e.
\begin{equation} \label{5.4f}
E  \left[ h_n(V_i,V_j) | V_i \right]=E  \left[ h_n(V_i,V_j) | V_j \right]=0~~~~~~ a.s.
\end{equation}
and apply Theorem 2.1 in de Jong (1987). For this purpose we calculate the variance as
\begin{eqnarray*}
\sigma_n^2 &=& \Var (\tilde{T}_n(c) ) = \frac{2^2}{n^2(n-1)^2} \sum_{i<j} E  [h_{ij}^2 ] \\
&=&  \frac{2}{n(n-1)}  E \left[ K_g^2(X_i-X_j) m^4(X_i)w^2(X_i)  m^4(X_j)  E[a^2(\varepsilon_i)|X_i] E[a^2(\varepsilon_j)|X_j ] w^2(X_j)\right]  \\
&=&   \frac{\mu_0^2}{n(n-1) g}  + o\left( \frac{1}{n^2 g}\right),
\end{eqnarray*}
where we have used the notation
\begin{equation}\label{ae}
a(\varepsilon_i) = 1-\varepsilon_i^2 - 2 c \varepsilon_i
\end{equation}
and $\mu_0^2$  is defined in (\ref{1.2a}). Observing that by this calculation $\sigma^2_{ij}:= E  [h_{ij}^2 ] = \frac{1}{g} \frac{\mu_0^2}{2}
\cdot (1 + o(1))$ we therefore have
\begin{eqnarray*}
\sigma^{-2}_n \max^n_{i=1} \Bigl( \frac{1}{n^2 (n-1)^2} \sum_{1 \leq j \leq n} \sigma_{ij}^2 \Bigr)= O \Bigl( \frac{1}{n} \Bigr),
\end{eqnarray*}
which proves the first assumption in de Jong's (1987) Theorem 2.1. In order to establish the second assumption we derive the decomposition [see
de Jong (1987)]
\begin{equation} \label{5.4g}
E[\tilde{T}_n^4(c) ]= G_{I} + 6 G_{II} + 12 G_{III} + 24 G_{IV} + 6 G_{V},
\end{equation}
where
\begin{eqnarray*}
G_I &=& \frac{2^4}{n^4(n-1)^4}\sum _{i < j} E[h_{ij}^4], \\
G_{II} &=& \frac{2^4}{n^4(n-1)^4}\sum _{i < j < k} E[h_{ij}^2h_{ik}^2] + E[ h_{ji}^2h_{jk}^2]  + E[h_{ki}^2h_{kj}^2],\\
G_{III} &=& \frac{2^4}{n^4(n-1)^4}\sum _{i < j < k} E[h_{ij}^2h_{ki}h_{kj}] + E[ h_{ik}^2h_{ji}h_{jk}]  + E[h_{kj}^2h_{ij}h_{ik}],\\
G_{IV} &=& \frac{2^4}{n^4(n-1)^4} \sum _{i < j < k<l} E[h_{ij}h_{ik}h_{lj}h_{lk}] + E[ h_{ij}h_{il}h_{kj}h_{kl}]  + E[h_{ik}h_{il}h_{jk}h_{jl}],\\
G_{V} &=& \frac{2^4}{n^4(n-1)^4} \sum _{i < j < k<l} E[h_{ij}^2h_{kl}^2] + E[ h_{ik}^2h_{jl}^2]  + E[h_{il}^2h_{jk}^2].
\end{eqnarray*}
It is easy to see that [recall the notation of  $a(\varepsilon_i)$ in (\ref{ae})]
\begin{eqnarray*}
E [h_{12}^4]
&=&  \frac{1}{g^4} \int  \int K_g^4 ( x_1-x_2) m^8 (x_1) m^8(x_2)  E[a^4(\varepsilon_1)|X_1=x_1] E[a^4(\varepsilon_2)|X_2=x_2]f(x_1) f(x_2)\\
&~&~~~~~~~~~~~~~~~~~~~ w^4(x_1)  w^4(x_2)~dx_1~ dx_2 \\
&=& O\left(\frac{1}{g^3}\right),\\
\end{eqnarray*}
which gives
$$
G_{I} =\frac{2^4}{n^4(n-1)^4} \sum _{i < j} E [h_{ij}^4] = O \left( \frac{1}{n^6 g^3} \right).
$$
The other terms are estimated similarly, i.e. $G_{II} =O ( \frac{1}{n^5 g^2} )$, $G_{III}=O( \frac{1}{n^5 g^2})$, $G_{IV}=O(\frac{1}{n^4 g})$
and it follows from (\ref{5.4g}) that
\begin{equation} \label{5.4h}
E [\tilde{T}_n^4(c)] = 6 G_{V} + o\left(\frac{1}{n^4 g^2} \right).
\end{equation}
On the other hand a straightforward calculation shows
\begin{eqnarray*}
\frac{1}{n^4(n-1)^4} \sum_{i<j} \sigma_{ij}^4 = O\left( \frac{1}{n^6 g^3} \right),~~\frac{1}{n^4(n-1)^4} \sum_{i<j<k} \sigma_{ij}^2
\sigma_{ik}^2 = O\left( \frac{1}{n^5 g^2} \right),
\end{eqnarray*}
and we obtain
\begin{eqnarray*}
\sigma^4_n &=&  2 G_{V} + o \left( \frac{1}{n^4g^2} \right),
\end{eqnarray*}
which proves the second assumption in Theorem 2.1 of de Jong (1987), that is
\begin{eqnarray*}
E [\tilde{T}_n^4(c)] - 3 \sigma^4_n= o\left( \frac{1}{n^4g^2} \right).
\end{eqnarray*}
Now this theorem yields asymptotic normality of $n \sqrt{g} ~\tilde{T}_n(c)´$ , i.e. $ n \sqrt{g} ~\tilde{T}_n(c)
\stackrel{\mathcal{D}}{\rightarrow} \mathcal{N}(0,\mu_0^2)$, and the assertion of Theorem 3.1(a) follows from (\ref{5.4d}) and (\ref{5.4e}).

\bigskip

{\bf Proof of Theorem 3.1(b).} We first note that under a fixed alternative the statistic $T_{3n}(c)$ is not negligible. More precisely,  we
show at the end of the proof
\begin{equation} \label{5.4i}
T_{3n}(c) - E [T_{3n}(c)] = 2 T_{5n}(c) + o_p\left( \frac{1}{\sqrt{n}} \right)
\end{equation}
where
\begin{eqnarray*}
E [T_{3n}(c)] = h^2 \kappa_2 E[\Delta_c(X) m(X) m''(X)f(X)w^2(X)]  + o( h^2).
\end{eqnarray*}
We obtain for the statistic $T_{6n}$ in (\ref{5.1})
\begin{eqnarray*}
 Var(T_{6n})
&=& \frac{2}{n^2 (n-1)^2} \sum_{i \ne j} E \left[ K_g^2( X_i - X_j) m^2(X_i) \sigma^2(X_i) w^2(X_i) ~m^2(X_j) \sigma^2(X_j) w^2(X_j) \right]\\
&=& \frac{2}{n (n-1)g} E[ m^4(x) \sigma^4(x)  f(x)w^4(x)]  \left( \int K^2(u)du  \right) + o\Bigl(\frac{1}{n^2 g} \Bigr),
\end{eqnarray*}
which gives
\begin{equation} \label{5.4k}
\sqrt{n} ~T_{6n} \stackrel{p}{\longrightarrow} 0
\end{equation}
(note that the expectation of $T_{6n}$ vanishes).  A similar calculation yields
\begin{eqnarray} \nonumber
T_{4n}(c) - E[T_{4n}(c)]
&=& \frac{2}{n} \sum_{i=1}  \Delta_c(X_i, \varepsilon_i) \Delta_c(X_i)f(X_i) w^2(X_i)  + o_p \left( \frac{1}{\sqrt{n}} \right), \\
&& \label{5.4l} \\
T_{5n}(c) &=&  \frac{1}{n} \sum_{i=1}^n \Delta_c(X_i) m(X_i) \sigma(X_i) f(X_i) w^2(X_i)  \varepsilon_i + o_p \Bigl( \frac{1}{\sqrt{n}} \Bigr)
\nonumber
\end{eqnarray}
and $E[T_{4n}(c)] =  E[\Delta_c^2(X) f(X)w^2(X) ] + O(g^2). $ Consequently we obtain from (\ref{5.1}), (\ref{5.4}), (\ref{5.4k})  and
(\ref{5.4l}) with
$$
E[T_n(c)]= E[\Delta_c^2(X) f(X)w^2(X) ] + h^2 B(c) + o(h^2)
$$
the following stochastic expansion under a fixed alternative
\begin{eqnarray}  \label{5.4m}
T_n(c) - E[T_n(c)] & =&  ( T_{4n}(c) - E[\Delta_c^2(X) f(X)w^2(X) ]) +  4T_{5n}(c)  + o_p
\Bigl({1\over \sqrt{n} } \Bigr) \\
&=& Z_n(c) + o_p\Bigl({1\over \sqrt{n} } \Bigr), \nonumber
  \end{eqnarray}
    where the random variable $Z_n(c)$ is defined by
  \begin{eqnarray}  \label{5.4n}
 Z_n(c) &=& {2\over n} \sum_{i=1}^n  \Bigl[  \Delta_c(X_i, \varepsilon_i)\Delta_c(X_i) f(X_i) w^2(X_i) - E[\Delta_c^2(X) f(X)w^2(X) ]  \nonumber \\
&&  ~~~~~~~~~~~  + 2  m(X_i) \sigma(X_i) \Delta_c(X_i) f(X_i) w^2(X_i) \varepsilon_i \Bigr],
 \end{eqnarray}
  and the assertion of Theorem 3.1(b) follows by a standard application of the central
 limit theorem verifying Ljapunoff's condition and observing that the
 dominating term on the right hand side of (\ref{5.4m}) has expectation $0$ and
 variance $\frac{\mu_1^2(c)}{n}$.

\bigskip

The proof of Theorem 3.1 is now completed showing the remaining estimates (\ref{5.4}) and (\ref{5.4i}).

\bigskip

{\bf Proof of the estimate (\ref{5.4}).} We consider exemplarily the case $j=1$ (the other case is treated by similar arguments) and  obtain
the  decomposition
\begin{equation} \label{5.5}
T_{1n} = 4 T_{1n}^{(1)} + 4 T_{1n}^{(2)} + T_{1n}^{(3)}
\end{equation}
with
\begin{eqnarray*}
T_{1n}^{(1)} &=& \frac{1}{n(n-1)} \sum_{i \ne j} K_g( X_i - X_j)  m(X_i) \tilde{\delta}_n(X_i)w(X_i)   m(X_j) \tilde{\delta}_n(X_j) w(X_j),\\
T_{1n}^{(2)} &=& \frac{1}{n(n-1)} \sum_{i \ne j} K_g( X_i - X_j)  m(X_i) \tilde{\delta}_n(X_i) w(X_i) \tilde{\delta}_n^2(X_j)w(X_j),  \\
T_{1n}^{(3)} &=& \frac{1}{n(n-1)} \sum_{i \ne j} K_g( X_i - X_j)  \tilde{\delta}_n^2(X_i)w(X_i)  \tilde{\delta}_n^2(X_j)w(X_j)
\end{eqnarray*}
[recall the definition of $\tilde{\delta}_n(x)$ in (\ref{5.5aa})]. The terms $T_{1n}^{(\ell )}$ are all treated similarly and we consider again
only the case $\ell=1$. With the notation
\begin{equation}\label{5.4aa}
d_k(x) = \frac{K_h(X_k - x)}{f(x)}~
\end{equation}
and Lemma 2 in Yao and Tong (2000) it follows
\begin{equation}\label{5.5a}
T_{1n}^{(1)}= \left( T_{1n}^{(1.1)}+ \kappa_2 T_{1n}^{(1.2)} + \frac{\kappa_2^2}{4} T_{1n}^{(1.3)} \right) \{1 + o_p(1)\}
\end{equation}
with
\begin{eqnarray*}
T_{1n}^{(1.1)} &=& \frac{1}{n^3(n-1)} \sum_{ i \ne j,  k,l }
K_g( X_i - X_j) d_k(X_i)  m(X_i)w(X_i)  d_l(X_j)  m(X_j)w(X_j)\sigma(X_k) \varepsilon_k \sigma(X_l) \varepsilon_l,\\
T_{1n}^{(1.2)} &=& \frac{h^2 }{n^2(n-1)} \sum_{i \ne j,  k } K_g( X_i - X_j) d_k(X_i) m(X_i)w(X_i) \sigma(X_k) \varepsilon_k   m(X_j) m''(X_j) w(X_j),\\
T_{1n}^{(1.3)} &=& \frac{h^4 }{n(n-1)} \sum_{ i \ne j} K_g( X_i - X_j)  m(X_i) m''(X_i)w(X_i)  m(X_j) m''(X_j)w(X_j).
\end{eqnarray*}
The terms $ T_{1n}^{(1.k)}$, $k=1,2,3$, can now be treated by calculating expectation and variance. For example,
\begin{eqnarray*}
&~&E [ T_{1n}^{(1.1)}]\\
 &=&  \frac{1}{n^3(n-1)} \sum_{ i \ne j,  k } E \left[K_g( X_i - X_j)  d_k(X_i)  m(X_i)w(X_i) d_k(X_j)  m(X_j) w(X_j)\sigma^2(X_k)  \right] + O  \Bigl( {1\over n^2h^2 }\Bigr)  \\
&=& {1\over n} \int K_g (x_i - x_j) K_h(x_k - x_i) K_h (x_k - x_j)    \\
&& ~~~~~~~~~~~~ \times m(x_i)w(x_i)  m(x_j)w(x_j) \sigma^2(x_k) f(x_k) dx_i dx_j dx_k + O  \Bigl( {1\over n^2h^2 }\Bigr)  \\
&=& \frac{1}{nh} \int m^2(x) \sigma^2(x)  f(x) w^2(x)~ dx \int K^2(w)~dw  + o\left( \frac{1}{nh} \right)  + O  \Bigl( {1\over n^2h^2 }\Bigr). \\
&=& O\left(\frac{1}{nh} \right)+ O\left(\frac{1}{n^2h^2} \right) = O\left(\frac{1}{nh} \right).
\end{eqnarray*}
For the calculation of the variance of $T_{1n}^{(1.1)}$ we introduce the notation $V_i=(X_i,\varepsilon_i)$,
$$
t_{1n}^{(1.1)}(V_i,V_j,V_k,V_l) =K_g( X_i - X_j) d_k(X_i)  m(X_i) w(X_i) \sigma(X_k) \varepsilon_k d_l(X_j)   m(X_j)w(X_j)  \sigma (X_l)
\varepsilon_l
$$
and obtain
\begin{eqnarray}
\Var  (T_{1n}^{(1.1)})  & \le &  E[ (T_{1n}^{(1.1)})^2 ]  \nonumber \\
&&=  \frac{1}{n^6(n-1)^2}  \sum_{i \ne j,  k,l }  \sum_{p \ne q, r,s } E \Bigl[ t_{1n}^{(1.1)} (V_i,V_j,V_k,V_l)~~t_{1n}^{(1.1)}
(V_p,V_q,V_r,V_s) \Bigr] \label{5.7}
= O \left(\frac{1}{n^2 h^2} \right), 
\end{eqnarray}
where we have used the fact, that the sum in (\ref{5.7}) is dominated by those expectations, for which two indices of $k,l,r,s$ coincide.
Markov's inequality and the assumptions on the bandwidth now yield $T_{1n}^{(1.1)} = o_p( n^{-1}g^{-1/2})$. The terms $T_{1n}^{(1.2)}$ and
$T_{1n}^{(1.3)}$ in (\ref{5.5a}) are treated similarly, which implies  $T_{1n}^{(1)} = o_p( n^{-1}g^{-1/2})$. A similar argument for the
statistics $T_{1n}^{(2)}$ and $T_{1n}^{(3)}$ yields
 assertion
(\ref{5.4}) for  $j=1$. The second case  $j=2$ is treated in the same way.

\bigskip

{\bf Proof of the estimate (\ref{5.4i}).} With the representation
$$
\delta_n(x)= \hat{m}^2(x) - m^2(x) = 2 m(x) \tilde{\delta}_n(x) + \tilde{\delta}_n^2(x)
$$
[recall the definition of $\tilde \delta(x)$ in (\ref{5.5aa})] we have
\begin{equation} \label{5.8}
T_{3n}(c)=2 T_{3n}^{(1)}(c) + T_{3n}^{(2)}(c)
\end{equation}
with
\begin{eqnarray*}
T_{3n}^{(1)}(c) &=& \frac{1}{n(n-1)} \sum_{i \ne j} K_g(X_i - X_j) m(X_i)\tilde{\delta}_n(X_i) w(X_i) \Delta_c(X_j, \varepsilon_j) w(X_j),\\
T_{3n}^{(2)}(c) &=& \frac{1}{n(n-1)} \sum_{i \ne j} K_g(X_i - X_j) \tilde{\delta}_n^2(X_i) w(X_i) \Delta_c(X_j, \varepsilon_j) w(X_j).
\end{eqnarray*}
A similar calculation as used in the proof of (\ref{5.4}) yields
\begin{equation}
T_{3n}^{(2)}(c) = O_p \left( \frac{1}{nh} \right)=o_p \left( \frac{1}{\sqrt{n}}\right). \label{5.9}
\end{equation}
For the first term in (\ref{5.8}) we introduce a further decomposition applying Lemma 2 in Yao and Tong (2000) and obtain
\begin{equation} \label{5.10}
T_{3n}^{(1)}(c) = \left(T_{3n}^{(1.1)}(c) + \frac{\kappa_2}{2} T_{3n}^{(1.2)}\right)~ \{ 1 + o_p(1)\}
\end{equation}
with [recall the definition of $d_k(\cdot)$ in (\ref{5.4aa})]
\begin{eqnarray*}
T_{3n}^{(1.1)}(c) &=& \frac{1}{n^2(n-1)} \sum_{i \ne j,  k } K_g(X_i - X_j) d_k(X_i) m(X_i) w(X_i) \sigma(X_k) \varepsilon_k \Delta_c(X_j, \varepsilon_j) w(X_j) ,\\
T_{3n}^{(1.2)}&=& \frac{h^2 }{n(n-1)} \sum_{i \ne j} K_g( X_i - X_j)  m(X_i) m''(X_i) w(X_i) \Delta_c(X_j, \varepsilon_j) w(X_j).
\end{eqnarray*}
 A straightforward but tedious calculation shows that $T_{3n}^{(1.1)}(c) $ and $T_{5n}(c)$ are
 asymptotically equivalent, that is
\begin{equation} \label{5.11}
 E [(T_{3n}^{(1.1)}(c) - T_{5n}(c))^2]=o\left(\frac{1}{n} \right).
\end{equation}
For the second term we obtain by a Taylor expansion
\begin{equation} \label{5.12}
 E [T_{3n}^{(1.2)}]= h^2  E[ \Delta_c(X) m(X)m''(X) f (X) w^2(X)] + o\left(\frac{1}{\sqrt{n}}\right)
\end{equation}
and
\begin{eqnarray*}
\Var  (T_{3n}^{(1.2)}) = o \left( \frac{1}{n} \right),
\end{eqnarray*}
which yields
\begin{equation} \label{5.13}
\sqrt{n} \left(T_{3n}^{(1.2)} - E [T_{3n}^{(1.2)}] \right)=o_p(1).
\end{equation}
Combining (\ref{5.8}) - (\ref{5.13}) and noting that $E[T_{3n}(c)] = \kappa_2 E [T_{3n}^{1.2} (c)] + o(h^2)$ now establishes (\ref{5.4i}) and
completes the proof of Theorem 3.1.  \hfill $\Box$

\subsection{Proof of Theorem 3.2.}

Note that $
\hat{c}^2 = \frac{R_1}{R_2}
$ where
$$
R_1 := \frac{1}{n}\sum_{i=1}^n \hat{m}^2(X_i) \hat{r}^2(X_i) w(X_i), \quad  R_2 := \frac{1}{n}\sum_{i=1}^n (\hat{\sigma}^2(X_i))^2 w(X_i).
$$
We now apply Lemma 2 in Yao and Tong (2000) and Lemma 1 in Fan and Yao (1998), which give asymptotic representations
 of $\hat m (x) - m(x)$  and $ \hat \sigma^2 (x) - \sigma^2 (x)$. A straightforward but tedious calculation yields
\begin{eqnarray} \label{6.2}
~~~~~~~~ R_1 - E[R_1] &=&  \frac{1}{n} \sum_{i=1}^n \left\{ \Bigl(m^2(X_i) \sigma^2(X_i) w(X_i) \varepsilon_i^2 - E[m^2(X) \sigma^2(X) w(X)]\Bigr) \right.\\
&~&~~~~~~~~~~~~~~~~~ \Bigl.+ ~ 2~  m(X_i) \sigma^3 (X_i) w(X_i) \varepsilon_i \Bigr\} + o_p\left( \frac{1}{\sqrt{n}} \right),\nonumber \\
~~~~~~~~  R_2 - E[R_2] &=& \frac{1}{n} \sum_{i=1}^n \biggl\{  \Bigl( \sigma^4(X_i) w(X_i) - E[\sigma^4(X) w(X)] \Bigr) \\
\nonumber && ~~~~~~~~~~~~~~~~~
 + 2 ~\sigma^4(X_i) w(X_i) \{ \varepsilon_i^2 -1 \}  \biggr\}+ o_p \left( \frac{1}{\sqrt{n}} \right)
\end{eqnarray}
with
\begin{eqnarray*}
E [R_1] &=& E[m^2(X) \sigma^2(X) w(X) ] + h^2 \kappa_2 E[ m(X) m''(X) \sigma^2(X) w(X)] + o ( h^2 ), \\
E [R_2] &= & E[\sigma^4(X) w(X)] + h^2 \kappa_2 E[\sigma^2(X) (\sigma^2(X))'' w(X)] + o (h^2).
\end{eqnarray*}
With this representation and the notations $r_1:=E[m^2(X) \sigma^2(X) w(X)]$ and $r_2:=E[\sigma^4(X) w(X)]$ it is easy to see that
 $$
 \frac{1}{R_2} = \frac{1}{r_2} \Bigl\{ 1 - \frac{R_2-r_2 }{r_2} + o \Bigl( \frac{1}{ \sqrt{n}} \Bigr) + o_p \Bigl( \frac{1}{ \sqrt{n}} \Bigr) \Bigr\},
$$
which implies
\begin{equation} \label{6.3}
\hat{c}^2- c_0^2 =\frac{R_1 - r_1}{r_2}-  \frac{R_1(R_2 -r_2 )}{ r_2^2} + o_p\Bigl( \frac{1}{ \sqrt{n}} \Bigr),
\end{equation}
where $c_0^2$ is defined in (\ref{2.6}). Observing (\ref{6.2}) we have
\begin{eqnarray*}
R_1(R_2 - r_2) &=&  E [R_1] (R_2 - r_2) + (R_1 - E [R_1]) (R_2  -r_2 ) \\
&{=}&E [R_1] (R_2  -E [R_2]) + E [R_1] (E [R_2] - r_2)+ O_p \Bigl( \frac{h^2}{\sqrt{n}}\Bigr)\\
&=& r_1 ( R_2 - E [R_2]) + E [R_1] (E [R_2]- r_2) + O_p \Bigl( \frac{h^2}{\sqrt{n}}\Bigr),
\end{eqnarray*}
and we obtain
\begin{eqnarray*}
\hat{c}^2- c_0^2& =& \frac{R_1 - E [R_1] }{r_2}  +  \frac{E[R_1] - r_1}{r_2} - \frac{r_1(R_2 - E [R_2])}{ r_2^2}
 -  \frac{E [R_1] (E [R_2] - r_2)}{ r_2^2}  + o_p\Bigl( \frac{1}{ \sqrt{n}} \Bigr) \\
 &=&h^2 \Gamma + o(h^2) + \tilde Z_n + o_p\Bigl( \frac{1}{ \sqrt{n}} \Bigr),
\end{eqnarray*}
where  the random variable $\tilde Z_n$ is defined as
\begin{eqnarray}
\tilde{Z}_n &=&  \frac{1}{n} \sum_{i=1}^n \Bigl\{ \tau_1 \Bigl( m^2 (X_i) \sigma^2(X_i) w(X_i) \varepsilon_i^2 - E[m^2(X) \sigma^2(X) w(X)]\Bigr)   + 2 \tau_1 m(X_i)  \sigma^3 (X_i) w(X_i) \varepsilon_i  \Bigr.  \nonumber\\
&~&~~~~~~~~~~~~~- \Bigl. \tau_2   \Bigl( \sigma^4(X_i) w(X_i) - E[\sigma^4(X) w(X)] \Bigr) -~2 \tau_2 ~ \sigma^4(X_i) w(X_i) \{ \varepsilon_i^2
-1 \} \Bigr\}.   \label{6.5}
\end{eqnarray}
The assertion finally follows from the central limit theorem and a straightforward but tedious calculation of the variance  $ n \Var
(\tilde{Z}_n)= \nu^2$. \hfill $\Box$

\subsection{Proof of Theorem 3.3.}

The assertion of the theorem follows, if the estimate
\begin{equation}\label{7.1}
T_n(c)  - T_n(\hat{c}) = o_p \left( \frac{1}{n \sqrt{g}} \right)
\end{equation}
can be established. For this purpose we introduce the analogue of the decomposition (\ref{5.1}), and obtain
\begin{eqnarray} \label{7.2}
~~~~~~
T_n (\hat{c}) &=&  (\hat{c}^2 + 1)^2 T_{1n} - 4 \hat{c}^2(\hat{c}^2 +1) T_{2n} + 2(\hat{c}^2 +1) T_{3n}(\hat{c}) + T_{4n}(\hat{c}) - 4 \hat{c}^2 T_{5n}(\hat{c}) + 4 \hat{c}^4 T_{6n} \\
&=& \nonumber 2(\hat{c}^2 +1) T_{3n}(\hat{c}) + T_{4n}(\hat{c}) - 4 \hat{c}^2 T_{5n}(\hat{c}) + 4 \hat{c}^4 T_{6n} + o_p\Bigl( {1\over n
\sqrt{g}} \Bigr),
\end{eqnarray}
where the last equality follows from (\ref{5.4}) and Theorem 3.2. For the third term in (\ref{7.2}) we have
\begin{eqnarray*}
T_{3n} (\hat{c}) = T_{3n}(c_0) - \left( \hat{c}^2 - c_0^2 \right)   T_{3n}^{(a)},
\end{eqnarray*}
where
\begin{eqnarray*}
T_{3n}^{(a)} = \frac{1}{n(n-1)} \sum_{i \ne j}  K_g( X_i - X_j ) \delta_n(X_i) w(X_i) \sigma^2(X_j)  w(X_j) \varepsilon_j^2 =  O_p(h^2) + O_p
\Bigl( \frac{1}{n h} \Bigr)
\end{eqnarray*}
[this estimate can be shown by similar arguments as presented  in the proof of the estimate (\ref{5.4})]. Consequently we obtain from  Theorem
3.2
\begin{equation}\label{7.3}
T_{3n}(\hat{c}) - T_{3n} (c_0) = o_p\Bigl(  \frac{1}{n \sqrt{g}}\Bigr).
\end{equation}
The corresponding estimates
\begin{eqnarray}\label{7.4a}
T_{4n}(\hat{c}) - T_{4n} (c_0) & \stackrel{H_0}{=}&  o_p\Bigl(  \frac{1}{n \sqrt{g}}\Bigr)~, \\
\hat{c}^2 T_{5n}(\hat{c}) &=& c_0^2 T_{5n} (c_0) + o_p \Bigl(  \frac{1}{n \sqrt{g}}\Bigr),
\label{7.4b} \\
(\hat{c}^2)^2 T_{6n}&=& c_0^4 T_{6n} + O_p \Bigl( \frac{h^2}{n \sqrt{g}} \Bigr) \label{7.4c}
\end{eqnarray}
are proved by similar arguments, which are  not given here for the sake of brevity. Note that the estimates (\ref{7.4b})
 and (\ref{7.4c}) hold also under a fixed alternative. Recalling
$$
c_0^2\stackrel{H_0}{=}c^2
$$
and combining  (\ref{7.2}) - (\ref{7.4b}) with Theorem 3.2 yields
\begin{equation}\label{7.5}
 T_n(\hat{c}) \stackrel{H_0}{=} T_{4n}({c}) - 4 {c}^2 T_{5n}({c}) + 4 {c}^4 T_{6n} + o_p\Bigl( {1\over n \sqrt{g}} \Bigr)
\stackrel{H_0}{=}  T_n({c}) + o_p\Bigl( {1\over n \sqrt{g}} \Bigr)
\end{equation}
where the last equality follows from the proof of Theorem 3.1(a). This proves the first identity in Theorem 3.3, while the remaining part
follows again from Theorem 3.1. \hfill $\Box$

\subsection{Proof of Theorem 3.4.}

Recalling the decomposition (\ref{7.2})  and (\ref{7.4b}), (\ref{7.4c}) we obtain
\begin{eqnarray} \label{7.8}
T_n (\hat{c}) &=&   2(\hat{c}^2 +1) T_{3n}(\hat{c}) + T_{4n}(\hat{c}) - 4 {c}_0^2 T_{5n}({c}_0) + 4 {c}_0^4 T_{6n}  + o_p\left(
\frac{1}{\sqrt{n}}\right),
\end{eqnarray}
and it remains to investigate the asymptotic properties of $T_{3n}(\hat{c})$ and $T_{4n}(\hat{c}) $ under the fixed alternative. For this we
note that (\ref{7.3}) and the obvious estimate $T_{3n}({c}_0) =o_p(1)$ yield
\begin{eqnarray} \label{7.9}
(\hat{c}^2 +1) T_{3n}(\hat{c}) - (c^2_0 +1) E[T_{3n}(c_0)] =2(c^2_0+1) T_{5n}(c_0)+  o_p\Bigl( \frac{1}{\sqrt{n}} \Bigr).
\end{eqnarray}
The corresponding estimate for the term $ T_{4n}(\hat{c})$ is more difficult. Note first that
\begin{equation} \label{7.10}
T_{4n}(\hat{c}) = T_{4n}(c_0) - 2 \left( \hat{c}^2 - c^2_0 \right) T_{4n}^{(a)}(c_0) + \left( \hat{c}^2 - c^2_0 \right)^2 T_{4n}^{(b)},
\end{equation}
where
\begin{eqnarray*}
T_{4n}^{(a)}(c_0) &\stackrel{H_1}{=}& \frac{1}{n(n-1)} \sum_{i \ne j}   K_g( X_i - X_j ) (m^2(X_i) - c_0^2 \sigma^2(X_i) \varepsilon_i^2) w(X_i) \sigma^2(X_j)  w(X_j) \varepsilon_j^2, \label{t4nunterh1t4na} \\
T_{4n}^{(b)} &\stackrel{H_1}{=}&  \frac{1}{n(n-1)} \sum_{i \ne j}   K_g( X_i - X_j ) \sigma^2(X_i)   w(X_i) \varepsilon_i^2 \sigma^2(X_j)
w(X_j) \varepsilon_j^2.  \label{t4nbh1}
\end{eqnarray*}
A standard calculation yields $T_{4n}^{(b)}=O_p(1)$, which shows that the second term in (\ref{7.10}) is of
 order $o_p({1\over \sqrt{n}})$. From Theorem 3.2. it follows
  $$
  (\hat{c}^2 - c_0^2) T_{4n}^{(a)} (c_0) = \varrho
  (\hat{c}^2 - E[\hat{c}^2])
  + h^2 \Gamma E[\Delta_{c_0}(X) \sigma^2(X) f(X) w^2(X)] + o( h^2)  + o_p\Bigl( \frac{1}{\sqrt{n}} \Bigr),
$$
where $\Delta_{c_0}(X)$ is defined in (\ref{2.4}). Combining this with (\ref{7.10}) therefore yields
\begin{eqnarray*}
T_n(\hat{c}) - E[T_n(\hat{c})]  & =& T_{4n}(c_0) - E[T_{4n}(c_0) ] + 4 T_{5n}(c_0) -2 \varrho
  (\hat{c}^2 - E[\hat{c}^2])   + o_p\Bigl( \frac{1}{\sqrt{n}}  \Bigr)\\
  &=&  Z_n(c_0)- 2 \varrho \tilde Z_n + o_p\Bigl( \frac{1}{\sqrt{n}}  \Bigr) ,
  \end{eqnarray*}
  where the random variables  $Z_n(c_0)$ and $\tilde Z_n$ are defined in (\ref{5.4n}) and (\ref{6.5}), respectively, and
  $$
  E[T_n(\hat{c})]= E[\Delta_{c_0}^2(X) f(X)w^2(X)] + h^2(B(c_0) - 2 \varrho \Gamma) + o(h^2).
  $$
The variances of $Z_n(c_0)$ and $\tilde Z_n$ have been determined in the proof of Theorem
  3.1 and 3.2 and a straightforward calculation gives $n \Cov(Z_n(c_0), \tilde{Z}_n)=\upsilon^2(c_0)$,
  where $\upsilon^2(c_0)$ is defined in Theorem 3.4. The assertion now follows from
  a standard application of the central limit theorem. \hfill $\Box$

\vskip 1cm

{\bf Acknowledgements.}

The authors would also like to thank an  anonymous referee for  constructive comments
on an earlier version of
his paper.
The work of the authors was supported by the Sonderforschungsbereich 475, Komplexit\"atsreduktion in multivariaten Datenstrukturen (Teilprojekt
B1) and an NIH grant award IR01GM072876:01A1.
The work of G. Wieczorek was additionally supported by a grant of the Ruhr-University Research
School.

\vskip 1cm

{\LARGE{ References}}

{\small{
\bigskip
R.J. Carroll, D. Ruppert (1988). {\itshape{Transformation and weighting in regression.}} Chapman and Hall, London.

\medskip
H. Dette (2002). A consistent test for heteroscedasticity in nonparametric regression based on the kernel method. {\itshape{J. Stat. Plann.
Inference}}, 103, (1-2), 311-329.

\medskip
H. Dette, A. Munk (2003). Some methodological aspects of validation of models in nonparametric regression. {\itshape{Statistica Neerlandica}},
57(2), 207--244.

\medskip
H. Dette, C. von Lieres und Wilkau (2003). On a test for constant volatility in continuous time financial models. {\itshape{Financ. Stoch.}} 7,
363--384.

\medskip
H. Dette, I. Spreckelsen (2004). Some comments on specification tests in nonparametric absolutely regular processes. {\itshape{J. Time Ser.
Anal.}} 25. No. 2. 159--172.

\medskip
H. Dette, B. Hetzler (2008). A simple test for the parametric form of the variance function in nonparametric regression. To appear in:
\emph{Ann.\  Inst.\  Stat.\ Math.}

\medskip
G.K. Eagleson, H.G. M\"uller (1997). Transformations for smooth regression models with multiplicative errors. {\itshape{J. R. Statisti. Soc.
B.}}, 59(1), 173--189.

\medskip
R.F. Engle (1982). Autoregressive heteroscedasticity of estimates of the variance of U.K. inflation. {\itshape{Econometrica}}, 50(4),
987--1007.

\medskip
J. Fan, I. Gijbels (1996). {\itshape{Local polynomial modelling and its applications.}} Chapman and Hall, London.

\medskip
J. Fan, Q. Yao (1998). Efficient estimation of conditional variance functions in stochastic regression. {\itshape{Biometrika}}, 85(3),
645--660.


\medskip
P. Gozalo, O. Linton (2000). Local nonlinear least squares: using parametric information in nonparametric regression. {\itshape{Journal of
Econometrics}}, 99, 63--106.

\medskip
P. de Jong (1987). A central limit theorem for generalized quadric forms. {\itshape{Probability Theory and Related Fields}}, 75, 261--277.

\medskip
P. Mc Cullagh, J. Nelder (1989). {\itshape{Generalized Linear Models}}. $2^{\scriptsize{\mbox{nd}}}$ ed. Chapman and Hall, London.

\medskip
A. Mil\o j (1985). The moment structure of ARCH processes. {\itshape{Scandinavian Journal of Statistics}}, 12, 281--292.

\medskip
G. Wieczorek (2007). Tests auf multiplikative Struktur in strikt station{\"a}ren Prozessen. {\itshape{PhD-Thesis, Ruhr University Bochum (in
German)}}.

\medskip
Q. Yao, H. Tong (2000). Nonparametric estimation of ratios of noise to signal in stochastic regression. {\itshape{Statistica Sinica}}, 10,
751--770.

\medskip
J.X. Zheng (1996). A consistent test of functional form via nonparametric estimation techniques. {\itshape{Journal of Econometrics}}, 75,
263--289. }}

\end{document}